\documentclass[MSNbibl,number,citesort,seceqn,dvips]{arxbj}
\usepackage{accents}

\aid{0}
\volume{18}
\issue{4}
\pubyear{2012}
\firstpage{1150}
\lastpage{1171}
\doi{10.3150/11-BEJ377}

\makeatletter
\newcommand{\R}{\mathbb{R}}
\newcommand{\N}{\mathbb{N}}
\newcommand{\F}{\mathcal{F}}
\newcommand{\M}{\accentset{\circ}{\mathcal{M}} }
\newcommand{\Nc}{\mathcal{N}}
\newcommand{\E}{\mathcal{E}}
\newcommand{\Es}{\mathbb{E}}
\newcommand{\Pro}{\mathbb{P}}
\newcommand{\di}{\,\mathrm{d}}
\large
\newcommand{\C}{\mathcal{C}}
\newcommand{\Ii}{\mathcal{I}}
\newcommand{\G}{\mathcal{G}}
\newcommand{\B}{\mathcal{B}}
\newcommand{\LL}{[\hspace{-0.5mm}[}
\newcommand{\comp}{\circ}
\newtheorem{theore}{Theorem}[section]
\newproclaim{defi}[theore]{Definition}
\newtheorem{propo}[theore]{Proposition}
\newremark{nota}[theore]{Remark}
\newtheorem{coro}[theore]{Corollary}
\newtheorem{lema}[theore]{Lemma}
\newcommand{\fraca}[2]{{#1}/{#2}}
\renewcommand{\backslash}{\setminus}
\makeatother

\begin{document}
\begin{frontmatter}

\title{Extended It\^o calculus for symmetric Markov processes}
\runtitle{Extended It\^o calculus for symmetric Markov processes}

\begin{aug}
\author{\fnms{Alexander} \snm{Walsh}\corref{}\ead[label=e1]{alexander.walsh@upmc.fr}}

\runauthor{A. Walsh}
\address{Laboratoire de Probabilit\'es et Mod\`eles Al\'eatoires, UMR
7599, Universit\'e
Paris VI, 4, Place Jussieu, 75252, Paris Cedex 05, France. \printead{e1}}
\end{aug}

\received{\smonth{7} \syear{2010}}
\revised{\smonth{5} \syear{2011}}

%
\begin{abstract}
Chen, Fitzsimmons, Kuwae and Zhang (\textit{Ann. Probab.} \textbf{36}
(2008) 931--970) have established an It\^{o} formula consisting in the
development of $F(u(X))$ for a symmetric Markov process $X$, a function
$u$ in the Dirichlet space of $X$ and any $\C^2$-function $F$. We give
here an extension of this formula for $u$ locally in the Dirichlet
space of $X$ and $F$ admitting a locally bounded Radon--Nikodym
derivative. This formula has some analogies with various extended It\^
{o} formulas for semi-martingales using the local time stochastic
calculus. But here the part of the local time is played by a process
$(\Gamma^a_t, a \in\R, t \geq0)$ defined thanks to Nakao's operator
(\textit{Z.~Wahrsch. Verw. Gebiete} \textbf{68} (1985) 557--578).
\end{abstract}

%
\begin{keyword}
\kwd{additive functional} \kwd{Fukushima decomposition} \kwd{It\^o formula}
\kwd{stochastic calculus} \kwd{symmetric Markov process} \kwd{zero
energy process}
\end{keyword}

\end{frontmatter}
%

\section{Introduction and main results}\label{sec1}\label{introduccion}\label{intc3}

For any real-valued semimartingale $Y = (Y_0+M_t+N_t)_{ t\geq0}$ ($M$
martingale and $N$ bounded variation process) and any function $F$ in
$\C^2(\R)$, the classical It\^{o} formula
\begin{eqnarray}
\label{itoc}
F(Y_t)&=&F(Y_0)+\int_0^t F^{\prime}(Y_s)\di M_s+\int_0^t F^{\prime
}(Y_s)\di N_s+\frac{1}{2}\int_0^t F^{\prime\prime}(Y_s)\di\langle
M^c \rangle_s\nonumber
\\[-8pt]
\\[-8pt]
&&{}+\sum_{s\leq t} \{F(Y_s)-F(Y_{s-})-F^{\prime}(Y_{s-})\Delta Y_s\}
\nonumber
\end{eqnarray}
provides both an explicit expansion of $(F(Y_t))_{t\geq0}$ and its
stochastic structure of semimartingale.

Let now $E$ be a locally compact separable metric space , $m$ a
positive Radon measure on $E$, and $X$ a $m$-symmetric Hunt process.
Under the assumption that the associated Dirichlet space $(\E,\F)$ of
$X$ is regular, Fukushima has showed that for any function $u$ in $\F
$, the additive functional (abbreviated as AF) $(u(X_t)-u(X_0))_{t \geq
0}$ admits the following unique decomposition:
\begin{equation}\label{fukushima}
u(X_t)=u(X_0)+M_t^u+N_t^u        \qquad   \Pro_x\mbox{-a.e. for
quasi-every $x$ in }   E,
\end{equation}
where $M^u$ is a martingale AF of finite energy and $N^u$ is a
continuous AF of zero energy.

Although $u(X)$ is not in general a semimartingale, Nakao~\cite{N} and
Chen \textit{et al}.~\cite{CFKZ} have proved that (\ref{itoc}) is still valid
with $u(X)$, $M^u$ and $N^u$ replacing, respectively, $Y$, $M$ and~$N$.
This is done thanks to the construction of a new stochastic integral
with respect to $N^u$, which takes the place of the well-defined
Lebesgue--Stieltjes integral for the bounded variation processes. As the
classical It\^{o} formula (\ref{itoc}), this It\^{o} formula for
symmetric Markov processes requires the use of $\C^2$-functions.

For the semimartingale case, there exist extended versions of (\ref
{itoc}) relaxing this regularity condition. This extensions are based
on the replacement of the fourth and fifth terms of the right-hand side
of (\ref{itoc}) by an alternative expression requiring only the
existence of $F^{\prime}$ and some integrability condition on
$F^{\prime}$ (see, e.g.,~\cite{E3,EK,EW}).
The integrability condition insures also the existence of the other
terms of (\ref{itoc}).

The question of relaxing the regularity condition on $F$ in the formula
of Nakao and Chen \textit{et al}. is a more complex question. Indeed the
integral   $\int_0^t F^{\prime}(u(X_s))\di N_s^u$   is well-defined
only when $F^{\prime}(u)$ belongs to $ \F_{\mathrm{loc}}$, the set of
functions locally in $\F$. As in~\cite{CFKZ}, $u\in\F_{\mathrm{loc}}$ means
that there exists a nest of finely open Borel sets $\{G_k\}_{k\in\N}$
and a sequence $\{u_k\}_{k\in\N}\subset\F$ such that $f=f_k$ q.e. on
$G_k$. As an example, in the case $X$ is a Brownian motion, this
condition implies that the second derivative $F^{\prime\prime}$
exists at least as a weak derivative.
Nevertheless, in the general case, we know that for any function $F$
element of $\C^1(\R)$ with bounded derivative, $F(u)$ belongs to $\F
$ and the process $F(u(X))$ hence admits a Fukushima decomposition. We
can thus hope to obtain an It\^{o} formula for $\C^1$-functions $F$
that would express each element of the decomposition of $F(u(X))$ in
terms of $F$, $u$, $N^u$ and $M^u$. 
Our purpose here is to establish such a formula. The obtained formula
is actually established for the functions $F$ with locally bounded
Radon Nikodym derivative and $u$ element of $\F_{\mathrm{loc}}$.

Before introducing this extended It\^{o} formula for symmetric Markov
processes, remark that one can easily obtain an extended It\^{o}
formula in case $u(X)$ is a semimartingale. Indeed, under the
assumption that $X$ has an infinite life time, we note (see (3.4) in
\cite{CFKZ}) that $u(X)$ is then a reversible semimartingale and that
one can hence make use of~\cite{E3} or~\cite{ERV} to develop
$F(u(X))$. But in general, $u(X)$ is not a semimartingale.

The extended It\^{o} formula for symmetric Markov processes presented
here is based on the construction for a fixed $t>0$, of a stochastic
integral of deterministic functions with respect to the process
$(\Gamma_t^a(u))_{a\in\R}$, defined as follows.

For $u$ in $\F$, let $M^{u,c}$ be the continuous part of $M^{u}$. For
any real $a$ and $t\geq0$, we set
\[
Z_t^a(u) = \int_0^t 1_{\{u(X_{s})\leq a\}}\di M_s^{u,c}
\]
and define $\Gamma^a$ by
\[
\Gamma^a(u) = (\Gamma_t^a(u))_{t \geq0} = (\Gamma(Z^a(u))_t)_{t
\geq0} = \Gamma(Z^a(u)),
\]
where $\Gamma$ is the operator on the space of martingale AF with
finite energy constructed by Nakao~\cite{N} (its definition is
recalled in Section~\ref{preli}). The process $(\Gamma^a_t(u))_{t
\geq0}$ is hence an additive functional with zero energy.

In Section~\ref{sec2}, we will see that the definition of $\Gamma^a (u)$ can be
extended to functions $u$ in $ \F_{\mathrm{loc}}$. In that case, the process
$M^{u,c}$ is a continuous martingale AF on $\LL0,\zeta\LL$ locally
of finite energy and the process $(\Gamma^a_t(u))_{t \geq0}$ is an AF
on $\LL0,\zeta\LL$ locally with zero energy.


As shown by the Tanaka formula (\ref{tanaka1}) below, the
doubly-indexed process $(\Gamma_t^a(u), a \in\R, t \geq0)$
plays almost the part of a local time process for $u(X)$. In Section~\ref{sec5},
this analogy with local time will be fully clarified under some
stronger assumption on $u$.



To introduce the obtained It\^{o} formula, we need the objects
presented by the following lemma. We denote by $(N(x,\mathrm{d} y),H)$ a L\'
{e}vy system for $X$ (See Definition A.3.7 of~\cite{FOT}), by $\nu_H$
the Revuz's measure of $H$ and by $\zeta$ the life time of $X$.

\begin{lema}\label{defZF}
Let $u\in\F$ (resp., $u\in\F_{\mathrm{loc}}$).
There exists a sequence $(\varepsilon_n)_{n\in\N}$ of positive real
numbers converging to $0$ and such that for any locally absolutely
continuous function $F$ from $\R$ into $\R$ with a locally bounded
Radon--Nikodym derivative, the following two processes are well-defined.
\begin{eqnarray}\label{firstlemma}
M_t^d(F,u)& = &\lim_{n\rightarrow\infty} \biggl\{\sum_{s\leq t} \{
F(u(X_s))-F(u(X_{s-}))\}1_{\{\varepsilon_n<|u(X_s)-u(X_{s-})|<1\}}1_{\{
s<\zeta\}} \nonumber\\
&&\hphantom{\lim_{n\rightarrow\infty} \biggl\{}{}- \int_0^t\int_{\{\varepsilon_n<|u(y)-u(X_s)|<1\}}\{
F(u(y))-F(u(X_s))\}N(X_s,\mathrm{d} y)\di H_s \biggr\}\nonumber\\
\label{f01120310}
A_t(F,u)& = &\lim_{n\rightarrow\infty}\int_0^t\int_{\{\varepsilon
_n<|u(y)-u(X_s)|<1\}}\{F(u(y))-F(u(X_s))\}N(X_s,\mathrm{d} y)\di H_s.\nonumber
\end{eqnarray}
The above limits are uniform on any compact of $[0,\infty)$ (resp.,
$[0,\zeta)$) $\Pro_x$-a.e. for q.e. $x\in E$. Moreover,
$(M_t^d(F,u))_{t\geq0}$ is a local martingale AF (resp., AF on $\LL
0,\zeta\LL$) with locally finite energy and $(A_t(F,u))_{t \geq0}$
is a continuous AF (resp., AF on $\LL0,\zeta\LL$) locally with $0$
energy.
\end{lema}


With the notation of Lemma~\ref{defZF}, we have the following It\^{o} formula.

\begin{theore}\label{itoSMP} Let $u\in\F$ (resp., $u\in\F_{\mathrm{loc}}$).
For any locally absolutely continuous function $F$ from $\R$ into $\R
$ with a locally bounded Radon--Nikodym derivative $F'$ such that $F(0)
= 0$,
the process $(F(u(X_t),t\in[0,\infty))$ (resp., $t\in[0,\zeta)$)
admits the following decomposition $\Pro_x$-a.e. for q.e. $x\in E$
\begin{equation}
F(u(X_t)) = F(u(X_0))+M_t(F,u)+Q_t(F,u)+V_t(F,u),
\label{fitoSMP}
\end{equation}
where $M(F,u)$ is a local martingale AF (resp., AF on $\LL0,\zeta\LL
$) locally of finite energy, $Q(F,u)$ is an AF (resp., AF on $\LL
0,\zeta\LL$) locally of zero energy,
and $V(F,u)$ is a bounded variation process, respectively, given by:
\begin{eqnarray}
M_t(F,u)&=&M_t^d(F,u)+\int_0^t F'(u(X_s))\di M^{u,c}_s,\nonumber\\
Q_t(F,u)&=&\int_{\R} F'(z)\di_{z} \Gamma_t^{z}(u)+A_t(F,u),\nonumber
\\
V_t(F,u)&=&\sum_{s\leq t}\{F(u(X_s))-F(u(X_{s-}))\}1_{\{
|u(X_s)-u(X_{s-})|\geq1\}}1_{\{s<\xi\}}\nonumber\\
&&{}-F(u(X_{\xi-}))1_{\{t\geq\xi\}}. \nonumber
\end{eqnarray}
\end{theore}

Note that for $u$ element of $\F$ and $F$ in $\C^2(\R)$, (\ref
{fitoSMP}) provides the It\^{o} formula of Chen \textit{et al}.~\cite{CFKZ}
together with the identity connecting integration with respect to
$(N^u_t)_{t \geq0}$ and integration with respect to $(\Gamma
^a_t(u))_{a \in\R}$ for smooth enough functions.

As a consequence of
Theorem~\ref{itoSMP}, we obtain the following Tanaka formula for
$\Gamma_t^a$:
\begin{eqnarray}\label{tanaka1}
\Gamma^a_t(u)&=&\bigl(u(X_0)-a\bigr)^{-}-\bigl(u(X_t)-a\bigr)^{-}+\int_0^t 1_{\{
u(X_{s-})\leq a\}}\di M_s^{u,c}\nonumber
\\[-8pt]
\\[-8pt] 
&&{}+\lim_{n\rightarrow\infty}\sum_{s\leq t}\bigl\{
\bigl(u(X_s)-a\bigr)^--\bigl(u(X_{s-})-a\bigr)^-\bigr\}1_{\{|u(X_s)-u(X_{s-})|>\varepsilon_n\}
},
\nonumber
\end{eqnarray}
where $(\varepsilon_n)_{n \in\N}$ is the sequence of Lemma \ref
{defZF} and the limit is uniform on any compact $\Pro_x$-a.e. for q.e.
$x\in E$. Using Tanaka's formula for semi-martingales (see~\cite{P}),
we obtain that when $u(X)$ is a martingale, $-2 \Gamma^a(u)$ is the
local time process of $u(X)$ at level $a$. This is the case when
$u(x)=x$ and $X$ is a symmetric L\'{e}vy process.

Formula (\ref{fitoSMP}) is hence reminiscent of various extensions of
It\^{o} formula involving stochastic integrals with respect to local
time, as for example the extensions given in~\cite{BY} for some
martingales,~\cite{E} for the Brownian Motion,~\cite{E2} and \cite
{EW} for L\'{e}vy processes with Brownian component and~\cite{W1} for
L\'{e}vy processes without Brownian component. Note that in case the
martingale part of $u(X)$ has no continuous component, the process
$\Gamma^a(u)$ is identically equal to $0$. But (\ref{fitoSMP}) still
represents an improvement of Fukushima's decomposition since (\ref
{fitoSMP}) requires only $u$ in $\F_{\mathrm{loc}}$ and $F$ with a locally
bounded Radon--Nikodym derivative.

Integration with respect to $(\Gamma^a_t(u))_{a \in\R}$ is
constructed in Section~\ref{sec3} and the It\^{o} formula (\ref{fitoSMP}) is
established in Section~\ref{sec4}.

In Section~\ref{SLT}, we will show that, when $\Gamma(M^{u,c})$ is of
bounded variation, $u(X)$ admits a local time process $(L^a_t, {a\in\R
}, t <\zeta)$ satisfying an occupation time formula of the same type
as the occupation time formula for the semimartingales and in this
case, the process of locally zero energy $Q(F,u)$ can be rewritten as:
\[
Q_t(F,u)=-\frac{1}{2}\int_{\R}F'(z)\di_z L_t^z+\int_0^t
F'(u(X_s))\di\Gamma(M^{u,c})_s+A_t(F,u), \qquad  t<\zeta.
\]

Finally in Section~\ref{SMC} we give a multidimensional version of
Theorem~\ref{itoSMP}.

\section{Preliminaries on $m$-symmetric Hunt
processes}\label{sec2}\label{preli}

Let $E$ be a locally compact separable metric space, $m$ a positive
Radon measure on $E$ such that $\operatorname{Supp}[m]=E$, $\Delta$ be a point
outside $E$ and $E_{\Delta}=E\cup\Delta$. Let $X=\{\Omega, \F
_{\infty},\F_t,X_t,\theta_t,\zeta,\Pro_x,x\in E_{\Delta}, t\geq
0\}$ be a $m$-symmetric Hunt Processes such that its associated
Dirichlet space $(\E,\F)$ is regular on $L^2(E;m)$. We may take as
$\Omega$ the space $D([0,\infty[\rightarrow E_{\Delta})$ of c\`
{a}dl\`{a}g functions from $[0,\infty[$ to $E_{\Delta}$, for which
$\Delta$ is a cemetery (i.e., if $\omega(t)=\Delta$, then $\omega
(s)=\Delta$ for any $s>t$) and denote by $\theta_t$ the operator
$\omega(s)\rightarrow\theta_t\omega(s):=\omega(t+s)$.
Every element $u$ of $\F$ admits a quasi-continuous $m$-version. In
the sequel, the functions in $\F$ are always represented by their
quasi-continuous $m$-versions. We use the term ``quasi everywhere'' or
``q.e.'' to mean ``except on an exceptional set.''

We say that a subset $\Xi$ of $ \Omega$ is a defining set of a
process $A=(A_t)_{t\geq0}$ with values in $[-\infty,\infty]$, if for
any $\omega\in\Xi,t,s \geq0\dvt\theta_t\Xi\subset\Xi$,
$A_0(\omega)=0$, $A_{.}(\omega)$ is c\`{a}dl\`{a}g and finite on
$[0,\zeta[$,
\[
A_{t+s}(\omega)=A_t(\omega)+A_s(\theta_t(\omega))
\]
and $A_t(\omega_{\Delta})=0$, where $\omega_{\Delta}$ is the
constant path equal to $\Delta$. A $(\F_t)$-adapted process is an
additive functional if it has a defining set $\Xi\in\F_{\infty}$
admitting an exceptional set, that is, $\Pro_x(\Xi)=1$ for q.e. $x\in E$.

An $(\F_t)$-adapted process is an additive functional on $\LL0,\zeta
\LL$ or a local additive functional if it satisfies all the conditions
to be an additive functional except that the additive property
$A_{t+s}(\omega)=A_t(\omega)+A_s(\theta_t(\omega))$ is required
only for $t+s<\zeta(\omega)$.

Let $\F_{\infty}^m$ (resp., $\F_t^m$) be the $\Pro_m$-completion of
$\sigma\{X_s, 0\leq s<\infty\}$ (resp., $\sigma\{X_s, 0 \leq s\leq
t\}$). An $(\F_t)$-adapted process is an additive functional admitting
$m$-null set if it has a defining set $\Xi\in\F_{\infty}^m$ such
that $\Pro_x(\Xi)=1$ for $m$-a.e. $x\in E$.

The abbreviations AF, PAF, CAF, PCAF and MAF stand respectively for
``additive functional,'' ``positive additive functional,'' ``continuous
additive functional,'' ``positive continuous additive functional'' and
``martingale additive functional,'' respectively. Let $\M$ and $\Nc_c$
denote, respectively, the space of MAF's of finite energy and the space
of continuous additive functionals of zero energy $N$ such that $\Es
_x(|N_t|)<\infty$ q.e. for each $t>0$. Moreover, $\M^c$ denotes the
subset of continuous elements of $\M$ and $\M^d$ denotes the subset
of purely discontinuous elements of $\M$.

For $u\in\F$, the elements $M^u$ and $N^u$ of the Fukushima's
decomposition (\ref{fukushima}) are elements of, respectively,
$ \M$ and $\Nc_c$. We denote by
$M^{u,c}$, $M^{u,j}$ and $M^{u,\kappa}$, respectively, the continuous
part, the jump part and the killing part of $M^u$ (see Section 5.3 of
\cite{FOT}). This three martingales are elements of $\M$.

Let $\Gamma$ the linear operator from $\M$ to $\Nc_c$ constructed by
Nakao~\cite{N} in the following way. It is shown in~\cite{N} that for
every $Z\in\M$, there is a unique $w\in\F$ such that
\[
\E(w,v)+(w,v)_m=\tfrac{1}{2}\mu_{\langle M^v+M^{v,\kappa},Z\rangle
}(E)      \qquad \mbox{for every }   v\in\F,
\]
where $(w,v)_m =\int_E w(x)v(x)m(\mathrm{d} x)$ and $\mu_{\langle
M^v+M^{v,\kappa},Z\rangle}$ is the smooth signed measure
corresponding to $\langle M^v+M^{v,\kappa},Z\rangle$ by the Revuz
correspondence. The process $\Gamma(Z)$ is then defined by:
\[
\Gamma_t(Z) = N_t^{w}-\int_0^t w(X_s)\di s.
\]
This operator satisfies: $\Gamma(M^u)=N^u$ for $u\in\F$. Thus $N^u$
admits the decomposition:
%
\begin{equation}\label{Nc}
N^u = {}^{c\!}N^{u} + {}^{j\!}N^{u} + {}^{\kappa\!}N^{u},
\end{equation}
where for $p\in\{c,j,\kappa\}\dvt{}^pN^{u}=\Gamma(M^{u,p})$.


For a Borel subset $B$ of $E\cup\{\Delta\}$, it is known that $\tau
_B =\inf\{t>0 \dvt X_t\notin B\}$ and
$\sigma_B =\inf\{t>0 \dvt X_t\in B\}$ are ($\F_t$)-stopping times.

An increasing sequence of Borel sets $\{G_k\}$ in $E$ is called a nest if
\[
\Pro_x \Bigl(\lim_{k\rightarrow\infty}\tau_{G_k}=\zeta
\Bigr)=1 \qquad \mbox{for q.e. } x\in E.
\]

%
%
Let $\mathcal{D}$ be a class of AF's. We say that an AF (resp., AF on
$\LL0,\zeta\LL$) is locally in $\mathcal{D}$ and write $A\in
\mathcal{D}_{\mathrm{loc}}$ (resp., $A\in\mathcal{D}_{f\mbox{-}\mathrm{loc}}$) if
there exists a sequence $\{A^n\}$ in $\mathcal{D}$ and an increasing
sequence of stopping times $T_n$ with $T_n\rightarrow\infty$ (resp., a
nest $\{G_n\}$ of finely open Borel sets) such that $\Pro_x$-a.e. for
q.e. $x\in E$, $A_t=A_t^n$ for $t<T_n$ (resp., $t<\tau_{G_n}$).

Let $\{A^n\}$ be a sequence in $\mathcal{D}$ such that for $k>n$,
$\Pro_x$-a.e. for q.e. $x\in E$, $A^k_t=A_t^n$ for $t<\tau_{G_n}$, then
it is clear that the process
\[
\everymath{\displaystyle }
A_t:= \left\{
\begin{array}{l@{ \qquad }l}
A^n_t&\mbox{for }t<\tau_{G_n},\\
0&\mbox{for }t\geq\zeta
\end{array}\right.
\]
is a well-defined element of $\mathcal{D}_{f\mbox{-}\mathrm{loc}}$.
A Borel function $f$ from $E$ into $\R$ is said to be locally in $\F$
(and denoted as $f\in\F_{\mathrm{loc}}$), if there is a nest of finely open
Borel sets $\{G_k\}$ and a sequence $\{f_k\}_{k\in\N}\subset\F$
such that $f=f_k$ q.e. on $G_k$. This is equivalent to (see Lemma~3.1(ii) in~\cite{CFKZ}) there is a nest of closed sets $\{D_k\}$ and
a sequence $\{f_k\}_{k\in\N}\subset\F_b$ such that $f=f_k$ q.e. on
$D_k$. For a such $f$,
\[
\everymath{\displaystyle }
M_t^{f,c}:= \left\{
\begin{array}{l@{ \qquad }l}
M_t^{f_k,c}&\mbox{for } t<\sigma_{E\backslash
G_k},\\0&\mbox{for } t\geq\lim_{k\rightarrow\infty}\sigma
_{E\backslash G_k}
\end{array}\right.
\]
is well defined and belongs to $\M_{f\mbox{-}\mathrm{loc}}$ because, for
$n>k$, $M_t^{f_n,c}=M_t^{f_k,c}\ \forall t\leq\sigma_{E\backslash
G_k}$ $\Pro_x$-a.e. for q.e. $x\in E$. Indeed, the last property is
shown in Lemma 5.3.1 in~\cite{FOT} for $\tau_{G_k}$ instead of
$\sigma_{E\backslash G_k}$, we conclude with the following observation:

For a CAF $A$, and a Borel set $G\subset E$, $\Pro_x$-a.e. for q.e.
$x\in E$:
%
\begin{equation}
A_t=0 \qquad \mbox{for }t<\tau_G \quad \Leftrightarrow \quad  A_t=0 \qquad \mbox{for }t<\sigma
_{E\backslash G}.\label{obvio}
\end{equation}

Every $f\in\F_{\mathrm{loc}}$ admits a quasi-continuous $m$-version, so we may
assume that all $f\in\F_{\mathrm{loc}}$ are quasi-continuous and we set
$f(\Delta)=0$.

We use the following notation for a locally bounded measurable function
$f$ and a $(\F_t)_{t \geq0}$-semimartingale $M$:
\[
(f * M)_t   =   \int_0^t f(X_{s-})\,\mathrm{d}M_s.
\]

We will use repeatedly the following fact (see Theorem 5.6.2 in~\cite{FOT}):

For any $F$ in $\mathcal{C}^1(\R^d)$ ($d$ is a positive integer) and
$u_1,\ldots,u_d$ in $\F_b$, the composite function $Fu=F(u_1,\ldots,u_d)$
belongs to $ \F_{\mathrm{loc}}$ and
\begin{equation}
M^{Fu,c}=\sum_{i=1}^d F_{x_i}(u)*M^{u_i,c}.\label{5.6.20FOT}
\end{equation}
Chen \textit{et al}.~\cite{CFKZ} have extended Nakao's definition of the
operator $\Gamma$ to the set of locally square-integrable MAF. We keep
using the letter $\Gamma$ for this extension without possible
confusion since thanks to Theorem 3.6 of~\cite{CFKZ} on the set $\M$,
both definitions given in~\cite{CFKZ} and~\cite{N} agree $\Pro
_m$-a.e. on $\LL0,\zeta\LL$. For a continuous locally
square-integrable MAF $M$, $\Gamma(M)$ is defined to be the following
CAF admitting $m$-null set on $\LL0,\zeta\LL$:
\begin{equation}
\Gamma_t(M)=-\tfrac{1}{2}(M_t+M_t\comp r_t)  \qquad \mbox{for
}t\in[0,\zeta
[,
\label{f02010310}
\end{equation}
where the operator $r_t$ is defined by
\[
r_t(\omega)(s) = \omega((t - s)-)1_{\{0\leq s\leq t\}}+\omega(0)1_{\{
s>t\}}   \qquad \mbox{for a path }    \omega\in\{t<\zeta\}
\]
and $r_t(\omega):=\omega_{\Delta}$ for a path $\omega\in\{t\geq
\zeta\}$.

The continuity of $\Gamma(M)$ $\Pro_m$-a.e. on $[0,\zeta[$ is a
consequence of Theorem 2.18 in~\cite{CFKZ}.

For $f$ a bounded element of $\F$ and $M$ in $\M$, Nakao has defined
the stochastic integral of $f(X)$ with respect to $\Gamma(M)$. We use
here the extension of this definition set by Chen \textit{et al}.~\cite{CFKZ}
for $f$ in $\F_{\mathrm{loc}}$ and $M$ continuous locally square-integrable MAF
as follows:
%
\begin{eqnarray}\label{integral}
f*\Gamma(M)_t  =   \int_0^tf(X_{s-})\di\Gamma_s(M)&:=&\Gamma
_t(f*M)-\frac{1}{2}\langle M^{f,c},M \rangle_t.%
\end{eqnarray}
It is a CAF admitting $m$-null set on $\LL0,\zeta\LL$.

When $M\in\M$ and $f\in\F_{\mathrm{loc}}$ the integral $f*\Gamma(M)_t$ can
be well defined $\Pro_x$-a.e. for q.e. $x \in E$. In particular, the
process $( f*\Gamma(M)_t)_{t \geq0}$ is a local CAF of $X$ (Lemma 4.6
of~\cite{CFKZ}).


The argument developed by Chen \textit{et al}. to write ``q.e. $x\in E$''
instead of ``$m$-a.e. $x\in E$'' in the proof of their Lemma 4.6 in \cite
{CFKZ}, is sufficient to establish Lemma~\ref{lema01280210} below.

\begin{lema}\label{lema01280210}
Let $A$ be an AF of $X$ (resp., AF on $\LL0,\zeta\LL$). Let $G$ be a
measurable subset of $ E_\Delta$ (resp., $G\subset E$) and $\Xi:=\{
\omega\in\Omega\dvt A_t\geq0, \forall t< \tau_{G}\}$, then $\Pro
_x(\Xi)=1$ for $m$-a.e. $x\in E$ if and only if $\Pro_x(\Xi)=1$ for
q.e. $x\in E$.
\end{lema}

\begin{lema}\label{lema01240310}
Let $\{D_n\}$ be a nest of closed sets and $\sigma:=\lim
_{n\rightarrow\infty}\sigma_{E\backslash D_n}$. Let $(M_n)_{n \in\N
}$ be a sequence of
$ \M^c$ such that for $n<k$, $\Pro_x$-a.e. for q.e. $x\in E$,
$M_t^n=M^k_t$ if $t<\sigma_{E\backslash D_n}$. Define a continuous
locally square-integrable MAF $M$ by:
\[
\everymath{\displaystyle }
M_t= \left\{
\begin{array}{l@{ \qquad }l}
M^n_t&\mbox{on } t<\sigma_{E\backslash D_n},\\
0&\mbox{on } t\geq\sigma.
\end{array}\right.
\]
Then
$\Gamma_t(M)$ can be well defined for all $t$ in $[0,\infty)$ $\Pro
_x$-a.e. for q.e. $x \in E$, by setting
\begin{equation}
\everymath{\displaystyle }
\Gamma_t(M)= \left\{
\begin{array}{l@{ \qquad }l}
\Gamma_t(M^n)&\mbox{on } t<\sigma_{E\backslash
D_n},\\0&\mbox{on } t\geq\sigma.
\end{array}\right.
 \label{defGloc}
\end{equation}
Moreover, $\Gamma(M)$ belongs to $\Nc_{c,f\mbox{-}\mathrm{loc}}$.
\end{lema}

For $f$ element of $\F_{\mathrm{loc}}$, (\ref{integral}) shows then that
$f*\Gamma(M)$ is a well defined CAF on $\LL0,\zeta\LL$.

%
\begin{pf*}{Proof of Lemma~\ref{lema01240310}}
A consequence of the $m$-symmetry assumption on $X$ is that the measure
$\Pro_m$, when restricted to $\{t<\zeta\}$ is invariant under $r_t$,
so we have $\Pro_m$-a.e. on $t<\zeta$:
$M_t\circ r_t=M^n_t\circ r_t$ if $t\leq\tau_{D_n}\circ r_t$, but
since $D_n$ is closed, for any $\omega\in\Omega$ and $t<\zeta
(\omega)$:
$t\leq\tau_{D_n}(\omega)\Leftrightarrow t\leq\tau_{D_n}(r_t\omega
)$. Hence, it follows from (\ref{f02010310}) that (\ref{defGloc})
hold $\Pro_m$-a.e. on $\LL0,\tau_{D_n} \LL$. This shows also, with
Lemma~\ref{lema01280210} that if $l>n$, $\Pro_x$-a.e. for q.e. $x\in
E$: $\Gamma_t(M^n)=\Gamma_t(M^l)$ for $t\leq\tau_{D_n}$ (and
consequently for $t\leq\sigma_{E\backslash D_n}$ by (\ref{obvio})).
Hence, the right-hand side of (\ref{defGloc}) is well defined as a CAF
belongs to $\Nc_{c,f\mbox{-}\mathrm{loc}}$.
\end{pf*}

\begin{nota}
Lemma~\ref{lema01240310} shows that for any $u\in\F_{\mathrm{loc}}$,
$^cN^u:=\Gamma(M^{u,c})$ is an element of $\Nc_{c,f\mbox
{-}\mathrm{loc}}$\label{notaNu}.
\end{nota}

The above Lemma~\ref{lema01280210} and Theorem 4.1 of~\cite{CFKZ}
lead to the following lemma.

\begin{lema}\label{lema01240210}
Let $M$ be an element of $\M$ such that $\Gamma(M)$ is of bounded
variation on each compact interval of $[0,\zeta[$. Then for every
element $f$ of $\F_{\mathrm{loc}}$, $\Pro_x$-a.e. q.e. for $x\in E$, on
$t<\zeta$, $\int_0^t f(X_s)\di\Gamma_s(M)$ coincides with the
Lebesgue--Stieljes integral of $f(X)$ with respect to $\Gamma(M)$.
\end{lema}

For the reader's convenience, we recall the following result which is
Theorem 5.2.1 of~\cite{FOT} and Theorem 3.2 of~\cite{N}, the last
statement can be seen directly from their proofs. By $e(M)$, we denote
the energy of $M$.

\begin{theore}
Let $\{M^n\dvt n\in\N\}$ be a $e$-Cauchy sequence of $\M$. There exists
a unique element $M$ of $\M$ such that $e(M^n-M)$ converges to zero.
The subsequence $n_k$ such that there exists $C\in\R_+$ such that for
every $k$ in $\N$: $e(M-M^{n_k})<C2^{-4k}$, satisfies: $\Pro_x$-a.e.
for q.e. $x\in E$, $M^{n_k}_t$ and $\Gamma_t(M^{n_k})$ converge
uniformly on any finite interval of $t$ to $M_t$ and $\Gamma_t(M)$,
respectively.
\label{NakaoFuku}
\end{theore}

\section{\texorpdfstring{Integration with respect to $\Gamma^z$}
{Integration with respect to Gamma z}}\label{sec3}\label{SiwrtG}

We fix a function $u$ of $\F_{\mathrm{loc}}$. Let $\{D_k\}_{k\in\N}$ be a
nest of closed sets and $(u_k)_{k\in\N}$ be a sequence of bounded
elements of $\F$ associated to $u$ such that $u=u_k$ q.e. on $D_k$. Let
$\sigma:=\lim_{n\rightarrow\infty}\sigma_{E\backslash D_n}$.
For any real number $a$, define $Z^a=Z^a(u)$ by
\[
\everymath{\displaystyle }
Z_t^a= \left\{
\begin{array}{l@{ \qquad }l}\int_0^t1_{\{u_k(X_{s-})\leq a\}}\di
M_s^{u_k,c}&\mbox{for } t\leq\sigma_{E\backslash D_k},\\\noalign{\vspace*{2pt}}
0&\mbox{for } t\geq\sigma.
\end{array}\right.
\]
$Z^a$ is a MAF on $\LL0,\zeta\LL$ locally of finite energy.
In particular, when $u$ belongs to $\F$, $Z^a$ is in $\M^c$ for any
real $a$. By Lemma~\ref{lema01240310}, $\Gamma(Z^a)$ is well-defined
and belongs to $\Nc_{c,f\mbox{-}\mathrm{loc}}$.

\begin{nota}
For $u$ element of $\F$, we can choose $D_k$ such that
%
\begin{equation}\label{correct}
\sigma=\lim_{k\rightarrow\infty}\sigma_{E\setminus D_k}=\infty
, \qquad \mbox{$\Pro_x$-a.e. for q.e. $x\in E$}.
\end{equation}
\end{nota}

Indeed, in this case, take $u_k:=(-k)\vee u\wedge k$ and $G_k:=\{
x\dvt|u(x)|<k\}$, then it follows from the strict continuity of $u$ that
$\lim_{k\rightarrow\infty}\sigma_{E\setminus G_k}=\infty$ $\Pro
_x$-a.e. for q.e. $x\in E$. Therefore, the nest of closed sets $\{F_k\}
_{k\in\N}$ built in the proof of Lemma 3.1(ii) in~\cite{CFKZ}
satisfies the property (\ref{correct}) and $u=u_k$ q.e. on $F_k$.
Choose then, $\{D_k\}=\{F_k\}$.

\begin{defi}
The process $(\Gamma^a_t, {a\in\R}, t \geq0)$ is defined by $\Gamma
^a_t=\Gamma^a_t(u) =\Gamma_t(Z^a)$. \label{defGx}
\end{defi}

Consider an elementary function $f$, that is, there exists two finite
sequences $(z_i)_{0\leq i\leq n}$ and $(f_i)_{0\leq i\leq n-1}$ of real
numbers such that:
\[
f(z)=\sum_{i=0}^{n-1}f_i1_{(z_i,z_{i+1}]}(z).
\]
For such a function integration with respect to $\Gamma_t=\{\Gamma
_t^z; z\in\R\}$ is defined to be the following CAF on $\LL0,\zeta
\LL$:
\begin{equation}
\int_{\R}f(z)\di_z\Gamma_t^z  = \sum_{i=0}^{n-1}f_i(\Gamma
_t^{z_{i+1}}-\Gamma_t^{z_i}).\label{intwrtG}
\end{equation}
Thanks to the linearity property of the operator $\Gamma$ we have for
any elementary function~$f$:
\[
\int_{\R} f(z)\di_z \Gamma_t^{z}=\Gamma_t \biggl(\int
_0^{\cdot}f(u(X_s))\di M_s^{u,c} \biggr).
\]

For any $k\in\N$, we define the norm $\|\cdot\|_k$ on the set of
measurable functions $f$ from $\R$ into $\R$ by
%
\begin{equation}
\|f\|_k= \biggl(\int_{E}f^2(u_k(x))\mu_{\langle M^{u_k,c}\rangle}(\mathrm{d}
x) \biggr)^{1/2}.\label{norma}
\end{equation}
Let $\Ii_k$ be the set of measurable functions from $\R$ into $\R$
such that $\|f\|_k<\infty$.

On $\Ii=\bigcap_{k\in\N} \Ii_k$, we define a distance $d$ by
setting:
\[
d(f,g) = [f-g],
\]
where
%
\begin{equation} \label{metric}
[f]=\sum_{k=1}^{\infty}2^{-k}(1\wedge\|f\|_k) .
\end{equation}
%

Note that $\Ii$ contains the measurable locally bounded functions and
that the set of elementary functions is dense in $(\Ii, d)$. Indeed,
by a monotone class argument, we can show that if $f$ is bounded, for
any $n\in\N$, there exists $f_n$ elementary such that $\sup_{k\leq
n}\|f-f_n\|_k\leq2^{-n}$. Hence,
\[
\sum_{n=1}^{\infty}[f-f_n]\leq\sum_{n=1}^{\infty} \Biggl( \sum
_{k=1}^n 2^{-k}(1\wedge\|f-f_n\|_k)+2^{-n}  \Biggr)< 2.
\]
Consequently it is sufficient to
show that the set of bounded functions is dense in $\Ii$. By dominated
convergence, $\lim_{n\rightarrow\infty}[f-(-n)\vee f \wedge n]=0$
for any $f\in\Ii$.

Let $f$ be an element of $\Ii$. The MAF $W^k$ defined by: $W^k_t=\int
_0^tf(u_k(X_s))\di M^{u^k,c}_s$, has finite energy since: $e(W^k)=\frac
{1}{2} \|f\|_k^2$. Hence,
\[
\everymath{\displaystyle }
fu*M^{u,c}_s:=\left \{
\begin{array}{l@{ \qquad }l}
fu_k*M^{u_k,c}_s&\mbox{for } t<\sigma_{E\backslash
D_k},\\0&\mbox{for } t\geq\sigma,
\end{array}\right.
\]
belongs to $\M^c_{f\mbox{-}\mathrm{loc}}$ ($\M^c_{\mathrm{loc}}$ if $u\in\F$) and by
Lemma~\ref{lema01240310}, $\Gamma(fu*M^{u,c})$ is well defined and is
an element of $\Nc_{c,f\mbox{-}\mathrm{loc}}$ ($\Nc_{c,\mathrm{loc}}$ if $u\in\F$).

\begin{theore}\label{int}
The application defined by (\ref{intwrtG}) on the set of elementary
functions can be extended to the set $\Ii$. This extension, denoted by
$\int f(z)\di_z\Gamma^z$, for $f$ in $\Ii$, satisfies:
\begin{enumerate}[(ii)]
\item[(i)] $\int f(z)\di_z\Gamma^z_t=\Gamma_t(fu* M^{u,c})\ \forall t\geq0$, $\Pro_x$-a.e. for q.e. $x\in E$.

\item[(ii)] Let $(f_n)_{n \in\N}$ be a sequence $ \Ii$. Assume
that: $[f_n-f]\rightarrow0$. Then there exists a subsequence
$(f_{n_k})_{k \in\N}$ such that $(\int f_{n_k}(z)\di_z\Gamma
^z_t)_{k \in\N}$ converges uniformly on any compact of $[0,\zeta)$
($[0,\infty)$ if $u\in\F$) to $\int f(z)\di_z\Gamma^z_t$ $\Pro
_x$-a.e. for q.e. $x\in E$.
\end{enumerate}
\end{theore}

\begin{pf}
Elementary functions are dense in $\Ii$ and (i) holds for elementary functions.
It is sufficient to prove that that if $[f_n-f]$ converge to zero,
there exists a subsequence $n_k$ such that for any $p\in\N$, $\Gamma
(f_{n_k}u*M^{u,c})$ converges to $\Gamma(fu*M^{u,c})$ uniformly on any
compact of $[0,\sigma_{E\backslash D_p}[$. Let $n_k$ be such that
$[f_{n_k}-f]<2^{-4k}$ and $p\in\N$, hence $\|f-f_{n_k}\|_p\leq
2^p2^{-4k}$ for any $k>p/4$ and it follows from Theorem \ref
{NakaoFuku} that $\Gamma(f_{n_k}u_p*M^{u_p,c})$ converges uniformly on
any compact to $\Gamma(fu_p*M^{u_p,c})$ $\Pro_x$-a.e. for q.e. $x\in
E$. But thanks to (\ref{defGloc}), $\Gamma(f_{n_k}u_p*M^{u_p,c})$ and
$\Gamma(fu_p*M^{u_p,c})$ agrees on $t<\sigma_{E\backslash D_p}$ with
$\Gamma(f_{n_k}u*M^{u,c})$ and $\Gamma(fu*M^{u,c})$, respectively,
$\Pro_x$-a.e. for q.e. \mbox{$x\in E$}.
\end{pf}

We finish this section with a characterization of the set $\Ii$ when
$u$ belongs to $\F$. Let $\E^{(c)}$ be the local part in the
Beurling--Deny decomposition for $\E$ (see Theorem~3.2.1 of~\cite
{FOT}). $\E^{(c)}$ has the local property, hence with the same
argument used to proof Theorems~5.2.1 and~5.2.3 of~\cite{BH}, there
exists a function $U$ in $L^1(\R,\mathrm{d} x)$ such that for any function
$F$ in $\mathcal{C}^1$ with bounded derivatives $f$:
\[
\E^{(c)}(F(u),F(u))=\frac{1}{2}\int_{\R}f^2(x)U(x)\di x.
\]

Then thanks to (\ref{5.6.20FOT}) and Lemma 3.2.3 of~\cite{FOT},
\[
\int_E f^2(u(x))\mu_{\langle M^{u,c}\rangle}(\mathrm{d} x)=\int_{\R
}f^2(x)U(x)\di x,
\]
hence it follows by a monotone class argument that for any measurable
positive function $f$ we have:
\begin{equation}\label{U}\int_E f(u(x))\mu_{\langle M^{u,c}\rangle
}(\mathrm{d} x)=\int_{\R}f(x)U(x)\di x.
\end{equation}

\begin{lema}\label{lemaUloc}
For $u$ element of $\F$, the set $\Ii$ coincides with the set
$L^1_{\mathrm{loc}}(\R,U(x)\di x)$, where the function $U$ is defined by (\ref{U}).
\end{lema}

\begin{pf}
For $k$ integer, the function $u_k$ is defined be $(-k)\vee u \wedge
k$. Associate $U_k$ to $u_k$ as $U$ is associated to $u$. We have then:
$\|f\|^2_k=\int_{\R}f^2(x)U_k(x)\di x$ for any measurable function~$f$. In order to proof Lemma~\ref{lemaUloc}, it is sufficient to prove
that: $U_k(x)=1_{[-k,k]}U(x)$ for a.e. $x$ in~$\R$.

Let $f$ be a continuous function with support in $[-k,k]$ and set
$F(x):=\int_0^xf(z)\di z$. We have hence: $F(u(x))=F(u_k(x))$ for any
$x$ in $E$ and therefore $f(u_k)*M^{u_k,c}=f(u)*M^{u,c}$, indeed thanks
to (\ref{5.6.20FOT}) both martingales coincides with $M^{Fu_k,c}$
($=M^{Fu,c}$).

We have therefore: $\int_E f^2(u_k(x))\mu_{\langle M^{u_k,c}\rangle
}(\mathrm{d} x)=\int_E f^2(u(x))\mu_{\langle M^{u,c}\rangle}(\mathrm{d} x)$. This
shows that
\[
\int_{\R}f^2(x)U_k(x)\di x=\int_{\R}f^2(x)U(x)\di x
\]
for any function $f$ continuous with compact support in $[-k,k]$, hence
$U_k(x)=U(x)$ for a.e. $x$ in $[-k,k]$.

Now if $g$ is a continuous positive function with support in $\R
\setminus[-k,k]$ then:
\[
\int_{\R}g(x)U_k(x)\di x=\int_E g(u_k(x)) \mu_{\langle
M^{u_k,c}\rangle}(\mathrm{d} x)=0
\]
therefore $U_k(x)=0$ for a.e. $x$ in $\R\setminus[-k,k]$. This
finishes the proof.
\end{pf}

\section{It\^{o} formula}\label{sec4}\label{SIF}

In this section, we first prove Lemma~\ref{defZF} and then Theorem
\ref{itoSMP}.\vadjust{\goodbreak}
\begin{pf*}{Proof of Lemma~\ref{defZF}} Let $u$ be an element
of $\F_{\mathrm{loc}}$, thanks to the proof of Lemma 3.1 of~\cite{CFKZ}, there
exists a nest of finely open
Borel sets $\{\G_k\}_{k\in\N}$ and a sequence $\{u_k\}_{k\in\N}$
in $\F$ such that $u(x)=u_k(x)$ for q.e. $x\in\G_k$ and $\|u_k\|
_{\infty}<k$. Let $\phi\in L^1(E;m)$ such that $0<\phi\leq1$ and
for any $k$ let
\[
h_k(x):=\Es_x \biggl(\int_0^{\sigma_{E\backslash\G_k}}\mathrm{e}^{-t}\phi
(X_t)\di t \biggr),
\]
$G_k:=\{x\in E\dvt h_k(x)>k^{-1}\}$ and $g_k(x):=1\wedge(kh_k(x))$. For
any $k$, $G_k\subset\G_k$, thus $u(x)=u_k(x)$ for q.e. $x\in G_k$.
Moreover, by the proof of Lemma 3.8 of~\cite{Ku}, $\{G_k\}_{k\in\N}$
is a nest and we have: $0\leq g_k\leq1$, $g_k(x)=1$ q.e. on $G_k$,
$g_k(x)=0$ on $E\backslash\G_k$. Since $h_k$ is quasi-continuous, we
can suppose that each $G_k$ is finely open (Theorem 4.6.1 of~\cite{FOT}).
For any $k\in\N$, we have:
\begin{eqnarray}
&&\int_{G_{k}}\int_{\{|u(x)-u(y)|<1\}}|u(x)-u(y)|^2N(x,\mathrm{d} y)\nu
_H(\mathrm{d} x)\nonumber\\
&& \quad =\int_{G_{k}}|g_{k}(x)|^2\int_{\{|u(x)-u(y)|<1\}
}|u(x)-u(y)|^2N(x,\mathrm{d} y)\nu_H(\mathrm{d} x)\nonumber\\
&& \quad \leq2\int_{G_{k}}\int_{\{|u(x)-u(y)|<1\}
}|g_{k}(x)-g_{k}(y)|^2|u(x)-u(y)|^2N(x,\mathrm{d} y)\nu_H(\mathrm{d} x)\nonumber\\
&& \qquad {}+2\int_{\G_{k}\times\G_{k} \cap\{|u(x)-u(y)|<1\}
}|g_{k}(y)|^2|u(x)-u(y)|^2N(x,\mathrm{d} y)\nu_H(\mathrm{d} x)\nonumber
\\
&& \quad \leq 2\int_{E\times E}|g_{k}(x)-g_{k}(y)|^2N(x,\mathrm{d} y)\nu_H(\mathrm{d}
x)\nonumber\\
&& \qquad {}+2\int_{E\times E}|u_{k}(x)-u_{k}(y)|^2N(x,\mathrm{d} y)\nu_H(\mathrm{d}
x)\nonumber\\
&& \quad \leq4\E(g_k,g_k)+4\E(u_k,u_k)<\infty.\nonumber
\end{eqnarray}

Therefore, if for any $\varepsilon>0$, we set:
\[
S_{\varepsilon}=\sum_{k=1}^{\infty}2^{-k} \biggl( 1 \wedge\int
_{G_k}\int_{\{|u(x)-u(y)|<\varepsilon\}}  |u(x)-u(y)|^2
N(x,\mathrm{d} y)\nu_H(\mathrm{d} x)  \biggr).
\]
We have then $\lim_{\varepsilon\rightarrow 0}S_{\varepsilon}=0$. We
choose a
sequence $(\varepsilon_n)_{n\in\N}$ such that $S_{\varepsilon_n}<2^{-4n}$.

Let $F$ be a locally absolutely continuous function with a locally
bounded Radon--Nikodym derivative $f$. For $k$ in $\N$, define $(F_k)$ by
\[
F_k(x)=F(x) 1_{ [-k-1,k+1]}(x) + F(k+1) 1_{[ k+1, \infty)}(x) +
F(-k-1) 1_{(-\infty, -k-1]}(x).
\]
Note that $F_k$ has a bounded Radon--Nikodym derivative: $f_k=f 1_{
[-k-1,k+1]}$.

For a function $\beta\dvtx E^2\rightarrow\R$, define:
\begin{eqnarray*}
A_t(\beta,n)&:=&\int_0^t\int_{\{\varepsilon_n<|u(y)-u(X_s)|<1\}}\beta
(y,X_{s})N(X_s,\mathrm{d} y)\di H_s \quad \mbox{and}
\\
M^d(\beta,n)&\hspace*{2.3pt}=&\sum_{s\leq t}\beta(X_s,X_{s-})1_{\{\varepsilon
_n<|u(X_{s-})-u(X_s)|<1\}}1_{\{s<\xi\}}-A_t(\beta,n).
\end{eqnarray*}
Denote by $M^d(F,u,n)$ (resp., $M^d(F,u,n,k)$) the process $M^d(\beta
,n)$ for $\beta(y,x)=F(u(y))-F(u(x))$ (resp., $\beta
(y,x)=(F(u(y))-F(u(x))1_{G_k}(x)$). Similarly, define $A^d(F,u,n)$ and
$A(F,n,u,k)$.

We just have to prove that $\Pro_x$-a.e. for q.e. $x\in E$, the limits
$\lim_{n\rightarrow\infty}M^d(F,u,n)$ and
$\lim_{n\rightarrow\infty}A(F,u,n)$ exist uniformly on any compact
of $[0,\sigma_{E\backslash G_k}[$. We have:
$M^d_t(F,u,n)=M^d_t(F_k,u,n,k)$ and $A_t(F,u,n)=A_t(F_k,u,n,k)$ on
$[0,\sigma_{E\backslash G_k}[$.
For every $k$, the process $M^d(F_k,u,n,k)$ belongs to $\M$ and for
$4n>k$, we have
\[
e\bigl(M^d(F_k,u,n+1,k)-M^d(F_k,u,n,k)\bigr)\leq c_k 2^k 2^{-4n},
\]
where $c_k=\|f_k\|_{\infty}$. Indeed, from the definition of
$\varepsilon_n$:
\begin{eqnarray}
&&e\bigl(M^d(F_k,u,n+1,k)-M^d(F_k,u,n,k)\bigr)\nonumber\\
&& \quad =\frac{1}{2}\int_{G_k\times E}\bigl(F_k(u(x))-F_k(u(y))\bigr)^21_{\{
\varepsilon_{n+1}\leq|u(x)-u(y)|<\varepsilon_n\}}N(x,\mathrm{d} y)\nu
_H(\mathrm{d} x)\nonumber\\
&& \quad \leq c_k\int_{G_k\times E}|u(x)-u(y)|^21_{\{|u(x)-u(y)|<\varepsilon
_n\}}N(x,\mathrm{d} y)\nu_H(\mathrm{d} x)\nonumber\\
&& \quad \leq c_k 2^k2^{-4n} \nonumber
\end{eqnarray}
thus, the convergence of $M^d(F,u,n)$ follows from Theorem \ref
{NakaoFuku}. Still thanks to Theorem~\ref{NakaoFuku}, the convergence
of $A(F,u,n)$ can be seen as a consequence of:
\begin{equation}
\Gamma(M^d_t(F_k,u,n,k))=A_t(F_k,u,n,k), \qquad      \Pro_x\mbox{-a.e. for
q.e. }  x\in E.
\label{f04150310}
\end{equation}
To prove (\ref{f04150310}), we note that $(A_t(F_k,u,n,k))_{t \geq0}$
is of bounded variation, so $A_t(F_k,u,n,k)\circ r_t=A_t(F_k,u,n,k)$
$\Pro_m$-a.e. on $t<\zeta$ (Theorem 2.1 of~\cite{Fi}). Hence, making
use of the operator $\Lambda$ defined in~\cite{CFKZ}, instead of
$\Gamma$, we first obtain:
\[
\Lambda(M^d_t(F_k,u,n,k))=A_t(F_k,u,n,k), \qquad      \Pro_m\mbox{-a.e. for
q.e. }  x\in E
    \mbox{ on }    \LL0,\zeta\LL.
\]
Finally by Theorem 3.6 in~\cite{CFKZ} and Lemma~\ref{lema01280210},
(\ref{f04150310}) holds, $\Pro_x$-a.e. for q.e. $x\in E$ on $\LL
0,\zeta\LL$, and therefore on $\LL0,\infty\LL$ thanks to the
continuity of $\Gamma(M^d_t(F_k,u,n,k))$ and $A_t(F_k,u,n,k)$.

It is clear that $M^d(F,u)\in\M_{f\mbox{-}\mathrm{loc}}$ and $A(F,u)\in\Nc
_{c,f\mbox{-}\mathrm{loc}}$. Moreover, for $u$ element of $\F$, we can take
$G_n=\{x\dvt|u(x)|<n\}$ for any $n$. In this case, from the strict
continuity of $u$ we have, $\Pro_x(\lim_{n\rightarrow\infty}\sigma
_{E\setminus G_n}=\infty)=1$ for q.e. $x\in E$, thus the convergence
of $M^{d}(F,u,n)$ and $A(F,u,n)$ are uniformly on any compact of
$[0,\infty)$. Thus, we obtain: $M^d(F,u)\in\M_{\mathrm{loc}}$ and $A(F,u)\in
\Nc_{c,\mathrm{loc}}$.
\end{pf*}

\begin{nota}\label{notajump}
(i) If $u\in\F$ and $f$ is bounded, then $M^d(F,u)\in\M$ and
$\Gamma(M^d(F,u))=A(F,u)$.
(ii)~With the notation of the proof of Lemma~\ref{defZF}, it holds
that if $u_k=u$ q.e. on~$G_k$:
\[
M^d_t(F,u)+A_t(F,u)=M^d_t(F_k,u_k)+A_t(F_k,u_k)    \quad\ \mbox{for }  t\in
[0,\sigma_{E\backslash G_k}[,\mbox{ $\Pro_x \mbox{-a.e. for q.e. }
x\in E$}.
\]
\end{nota}
%


%
\begin{pf*}{Proof of Theorem~\ref{itoSMP}}
We use the notation of the proof of Lemma~\ref{defZF}. Thus, if $u\in
\F$, we take $G_n:=\{x\dvt|u(x)|<n\}$, $n\in\N$. Let $F$ be a locally
absolutely continuous function $F$ with a locally bounded
Radon--Nikodym derivative $f$.

Let $I_t$ be the difference of the left-hand side and the right-hand
side of (\ref{fitoSMP}). For any $k$, we define $I_t^k$ as $I_t$ with
$u_k$ and $f_k$ replacing $u$ and $f$, respectively. Hence, $I_t=I_t^k$
for $t<\sigma_{E\backslash G_k}$, $\Pro_x$-a.e. for q.e. $x\in E$.
Since $\sigma_{E\backslash G_n}\wedge\zeta$ converges to $\zeta$ if
$u\in\F_{\mathrm{loc}}$ and $\sigma_{E\setminus G_n}$ converges to $\infty$
if $u\in\F$, it is sufficient to prove (\ref{fitoSMP}) on $[0,\sigma
_{E\backslash G_k}[$ for any $k\in\N$. Consequently, we can assume
(and we do) that $u$ is an element of $\F_b$ and $f$ is bounded.

If $f$ is continuous, thanks to (\ref{5.6.20FOT}), $F(u)\in\F$ and
$M^{Fu,c}=fu*M^{u,c}$ and we have the Fukushima decomposition:
\[
F(u(X_t))=F(u(X_0))+fu*M^{u,c}_t+\Gamma(fu*M^{u,c})_t+M^{u,d}_t+\Gamma
(M^{u,d})_t.
\]
We obtain (\ref{fitoSMP}) from Lemma~\ref{int}(i) and Remark \ref
{notajump}(i).

If $f$ is not necessarily continuous, let $g$ be in $L^1(\R)$ be a
strictly positive function on $\R$ such that $g$ and $1/g$ are locally
bounded. Define the norms $\|\cdot\|$ and $\|\cdot\|_*$ on the Borel
measurable functions as follows:
\begin{eqnarray*}
\|h\|_*&=& \biggl(\int_{E}h^2(u(x))\mu_{\langle M^{u,c}\rangle}(\mathrm{d}
x) \biggr)^{1/2},
\\
\|h\|&=&\|h\|_*+\int|h(x)|g(x)\di x\nonumber\\
&&{} + \biggl( \int_{E\times E-\delta}|u(x)-u(y)|\int
_{u(x)\wedge u(y)}^{u(x)\vee u(y)}h(z)^2\di z N(x,\mathrm{d} y)\nu_H(\mathrm{d} x)
 \biggr)^{\fraca{1}{2}}.\nonumber
\end{eqnarray*}
Since $u$ is in $\F$, we have $\|f\|<\infty$. By a monotone class
argument, one shows that there exists a sequence of bounded continuous
functions $(f_n)_{n\in\N}$ with compact support such that $\|f_n-f\|$
converges to $0$ as $n$ tends to infinity. We set $F_n(x) =\int_0^x
f_n(z)\di z$.

In order to show (\ref{fitoSMP}), we will show that there exists a
subsequence $n_k$ such that the terms in the expansion (\ref{fitoSMP})
for $F_{n_k}$ converge as $k\rightarrow\infty$ to the corresponding
expression with $f$ replacing $f_{n_k}$. The convergence of
$F_n(u(X_t))-F_{n}(u(X_0))-V_t(F_n,u)$ to
$F(u(X_t))-F(u(X_0))-V_t(F,u)$ is a consequence of the pointwise
convergence of $F_n$ to~$F$, indeed, for any $x\in\R$,
\[
|F_n(x)-F(x)|\leq\int_{-x^{-}}^{x^+}|f_n(z)-f(z)|\di z\leq\sup
_{|\lambda|\leq|x|}\frac{1}{g(\lambda)}\int_{-\infty}^{\infty
}|f_n(z)-f(z)|g(z)\di z\rightarrow0.
\]
The existence of a subsequence $\{n_k\}$ such that $\int_0^t
f_{n_k}(u(X_s))\di M_s^{u,c}$ and $\int_{\R}f_{n_k}(z)\di_z\Gamma
_t^z(u)$ converge to $\int_0^t f(u(X_s))\di M_s^{u,c}$ and $\int_{\R
}f(z)\di_z\Gamma_t^z(u)$, respectively, is a consequence of the fact
that $e(fu*M^{u,c}-f_nu*M^{u,c})=\frac{1}{2}\|f-f_n\|_*\rightarrow0$
as $n\rightarrow\infty$, and Theorem~\ref{NakaoFuku}. Thanks to
Theorem~\ref{NakaoFuku} and Remark~\ref{notajump}(i), it is then
sufficient to show that $e(M(F_n,u)-M(F,u))$ converges to zero as
$n\rightarrow\infty$. But
\begin{eqnarray}
e(M-M^n)&\leq&\frac{1}{2}\int_{E\times E-\delta
}\bigl(F(u(x))-F_n(u(x))-F(u(y))\nonumber\\
&&\hphantom{\frac{1}{2}\int_{E\times E-\delta
}\bigl(}{} +F_n(u(y))\bigr)^2N(x,\mathrm{d} y)\nu_H(\mathrm{d} x)\nonumber\\
&\leq& \frac{1}{2}\|f-f_n\|_*^2\rightarrow0 \qquad \mbox{as }n\rightarrow
\infty.\nonumber
\end{eqnarray}
\upqed
\end{pf*}
As an example, for $F(z)=z$ and $u$ in $\F_{\mathrm{loc}}$, one obtains a
Fukushima decomposition for the process $u(X)$. This case can be seen
as a refinement of Lemma 2.2 in~\cite{CFTY}.

\section{Local time}\label{sec5} \label{SLT}

We fix an element $u$ of $\F_{\mathrm{loc}}$. The associated process $^cN^u$
has been defined in (\ref{Nc}) by $^cN^{u}=\Gamma(M^{u,c})$. By
Remark~\ref{notaNu}, $^cN^{u}$ is a CAF locally of zero energy or
merely a CAF of zero energy when $u$ belongs to $F$. We suppose that
$u$ satisfies the additional assumption that $^cN^u$ is of bounded
variation on $[0,\zeta)$, that is, there exists two PCAF's $A^{(1)}$
and $A^{(2)}$ such that $\Pro_x$-a.e. for q.e. $x\in E$:
\begin{equation}
^cN_t^{u} =A_t^{(1)}-A_t^{(2)}  \qquad \forall t\in[0,\zeta
).\label{sup01}
\end{equation}
%

We remind that a measure $\nu$ on $E$ is a smooth signed measure on
$E$ if there exists a nest $\{F_k\}$ such that for each $k$,
$1_{F_k}.\nu$ is a finite signed Borel measure charging no set of zero
capacity and further $\nu$ charges no Borel subset of $E\setminus
\bigcup_{k=1}^{\infty}F_k$. Such nest is said to be associated to
$\nu$. For a closed set $F\subset E$, we set:
\[
\F_{b,F}=\{u\in\F_b\dvt u=0\mbox{ q.e. on }E\setminus F\}.
\]
We also need the following definition:
\[
\E_1(u,v)=\E(u,v)+(u,v)_m.
\]

\begin{lema}\label{Beurling} The process
$^cN^{u}$ is of bounded variation if and only if there exists a smooth
signed measure $\nu$ on $E$ with associated nest $\{F_k\}$ such that
\[
\E^{(c)}(u,v)=\langle v,\nu\rangle , \qquad    \forall v\in\bigcup
_{k=1}^{\infty}\F_{b,F_k}.
\]
\end{lema}

\begin{pf}
From Theorem 5.2.4 of~\cite{FOT}, $^cN^{u}$ is the only AF of zero
energy such that for any $h\in\F$,
\[
\lim_{t\downarrow0}\frac{1}{t}\Es
_{h.m}[^cN_t^{u}]=-e(M^{u,c},M^{h,c})=-\E^{(c)}(u,h).
\]
On the other hand, since: $|\E^{(c)}(u,h)|\leq(\E
^{(c)}(u,u))^{1/2}(\E_1(h,h))^{1/2}$, there exists a unique $w\in\F$
such that
\[
\E^{(c)}(u,h)=\E_1(w,h) \qquad \mbox{for any $h\in\F$}.
\]
Hence, $\lim_{t\downarrow0}\frac{1}{t}\Es_{h.m}[N_t^w-\int_0^t
w(X_s)\di s]=-\E^{(c)}(u,h)$ for all $h\in\F$. This implies that the
AF $N^w-\int_0^{\cdot} w(X_s)\di s$   is equivalent to $^cN^{u}$.
Consequently, $^cN^{u}$ is of bounded variation if and only if $N^w$ is
of bounded variation. But thanks to Theorem 5.4.2 of~\cite{FOT}, this
last condition is equivalent to the existence of a smooth signed
measure $\nu$ with an associated nest $\{F_k\}$ such that
\[
\E_1(w,v)=\langle v,\nu\rangle   \qquad \forall v\in\bigcup_{k=1}^{\infty}\F_{b,F_k}.
\]
\upqed
\end{pf}

\subsection{Definition of local time}\label{secu5.1}

\begin{defi}
The local time at $a$ of $u(X)$, denoted by $L_t^a=L_t^a(u)$ is the
following CAF on $\LL0,\zeta\LL$:
\[
\frac{1}{2}L_t^a =-\Gamma(Z^a)_t+\int_0^t 1_{\{u(X_{s-})\leq a\}}\di
^cN_s^{u}     \qquad \mbox{for $t\in[0,\zeta)$}.
\]
\end{defi}

The name ``local time'' is justified by Proposition~\ref{occuptime} and
Corollary~\ref{loctime} below.

\begin{propo}
\label{occuptime}
There exists a $\B(\R)\otimes\B(\R_{+})\otimes\F^m_{\infty
}$-measurable version of the local time process $\{\tilde{L}_t^a; a\in
\R,t\geq0\}$ such that $\Pro_m$-a.e. we have the occupation time
density formula:
\[
\int_{\R}f(x)\tilde{L}_t^{x}\di x=\int_0^tf(u(X_s))\di\langle
M^{u,c}\rangle_s     \qquad \mbox{for any $f$ Borel bounded and
$t<\zeta$}.
\]
\end{propo}

\begin{pf}
We start with the case when $u$ is a bounded element of $\F$. From
(\ref{f02010310}) we have: $\Pro_m$-a.e. on $\LL0,\zeta\LL$:
$L_t^a=Z_t^a+Z_t^a\comp r_t+2\int_0^t1_{\{u(X_{s-})\leq a\}}\di
N_s^{u,c}$. Moreover, thanks to Theorem 63, Chapter~IV of~\cite{P},
there exists a function $\tilde{Z}(a,t,\omega)$ in $\B(\R)\otimes
\B(\R_{+})\otimes\F^m_{\infty}$, such that for each $a\in\R$,
$\tilde{Z}(a,t,w)$ is a continuous $(\F^m_t)$-adapted version of the
stochastic integral $Z^a$, and thanks to Lemma 2.10 and Theorem 2.18 of
\cite{CFKZ}, $\tilde{Z}(a,t,w)\comp r_t(\omega)\in\B(\R)\otimes
\B(\R_{+})\otimes\F^m_{\infty}$ is a continuous $(\F
^m_t)$-adapted version of $Z^a_t\comp r_t$ for each $a\in\R$.
Besides, we can take $\int_0^t 1_{\{u(X_{s-})\leq a\}}\di N_s^{u,c}$
jointly continuous in $t$ and right continuous in $a$, $\Pro_m$-a.e. on
$\LL0,\zeta\LL\times\R$. Thus, we have constructed a version $\{
\tilde{L}_t^a, a\in\R, t\in[0,\zeta[\}$ of $\{L^a_t, a\in\R,t\in
[0,\zeta)\}$ which is $\B(\R)\otimes\B(\R_{+})\otimes\F
^m_{\infty}$-measurable.

Let $f$ be a continuous positive element of $L^1(\R)$. Using the proof
presented in~\cite{P} of Fubini's theorem for stochastic integrals
(Theorem 64, Chapter IV of~\cite{P}), we know that $\int_{\R}\tilde
{Z}(z,t,\omega)f(z)\di z$ is a well-defined Lebesgue integral since
$\Pro_m$-a.e.:
\[
\int_{\R}|\tilde{Z}(z,t,\omega)|f(z)\di z<\infty \qquad \mbox{for all }t.
\]
Moreover, still thanks to this theorem,
$\int_{\R}\tilde{Z}(z,t,\omega)f(z)\di z$ is a continuous $\Pro
_m$-version of
$\int_0^t F(u(X_s))\di M_s^{u,c}$, where $F(z)=\int_z^{\infty
}f(\lambda)\di\lambda$. Consequently, for $t>0$, $\Pro_m$-a.e. on $\{
t<\zeta\}$, $\int_{\R}|\tilde{Z}(z,t,r_t(\omega))| f(z)\di
z<\infty$ and $\int_{\R}\tilde{Z}(z,t,r_t(\omega)) f(z)\di z$ is a
continuous $\Pro_m$-version of
$\int_0^t F(u(X_s))\di M_s^{u,c}\comp r_t$.

Since $(\int_0^t 1_{\{u(X_s-)\leq a\}}\di N_s^{u,c})_{a \in\R}$ is
of bounded variation on $\{t<\zeta\}$, we obtain $\Pro_m$-a.e. on $\{
t<\zeta\}$: $ \int_{\R}f(z)|\tilde{L}_t^z|\di z < \infty$ and
\[
\int_{\R}f(z)\tilde{L}_t^z\di z=\int_0^t F(u(X_s))\di
M_s^{u,c}+\int_0^tF(u(X_s))\di M_s^{u,c}\comp r_t+2\int_0^t
F(u(X_{s}))\di N_s^{u,c}
\]
which leads to
%
\begin{equation}
\int_{\R}f(z)\tilde{L}_t^z\di z =-2\Gamma(Fu*M^{u,c})_t+2\int_0^t
F(u(X_{s}))\di N_s^{u,c}.\label{f2170210}
\end{equation}
Now thanks to (\ref{5.6.20FOT}), $Fu$ belongs to $\F_{\mathrm{loc}}$ and
$M^{Fu,c}_t=-\int_0^tf(u(X_s))\di M_s^{u,c}$. Thus,
\[
\langle M^{Fu,c},M^{u,c}\rangle_t=-\int_0^tf(u(X_s))\di\langle
M^{u,c}\rangle_s.
\]
Thanks to Lemma {\ref{lema01240210}} we have $\Pro_m$-a.e. on $\{
t<\zeta\}$:
\[
\int_0^tF(u(X_s))\di^cN_s^{u}=\int_0^tF(u(X_s))\di\Gamma(M^{u,c})_s.
\]
On the other hand, the definition of the integral with respect to
$\Gamma(M^{u,c})$ (Chen \textit{et al}.~\cite{CFKZ}) gives:
\[
\int_0^tF(u(X_s))\di\Gamma(M^{u,c})_s =\Gamma(Fu*M^{u,c})_t+\frac
{1}{2}\int_0^tf(u(X_s))\di\langle M^{u,c}\rangle_s
\]
which together with (\ref{f2170210}) lead to
\begin{equation}
\int_{\R}f(z)\tilde{L}_t^{z}\di z=\int_0^tf(u(X_s))\di\langle
M^{u,c}\rangle_s
, \qquad \Pro_m  \mbox{-a.e. on}  \{t<\zeta\}.
\label{nada}
\end{equation}
Actually, the set of null $\Pro_m$-measure on which (\ref{nada})
could fail can be chosen independently of~$f$. Indeed, the set of
continuous functions with compact support, is a separable topological
space for the metric of uniform convergence.

We show now that the set of null $\Pro_m$-measure on which (\ref
{nada}) could fail does not depend on $t$ either. We have thanks to
(\ref{nada})
\begin{equation}\Pro_m\mbox{-a.e. on }\{t<\zeta\}, \qquad  \tilde
{L}_t^z\geq0\mbox{ for }\mathrm{d} z\mbox{-a.e. }z\label{Zpos}
\end{equation}
hence by a monotone class argument, (\ref{nada}) holds $\Pro_m$-a.e.
on $\{t<\zeta\}$ for any $f$ Borel bounded. It remains to show that
(\ref{nada}) holds $\Pro_m$-a.e. on $\LL0,\zeta\LL$. To do so, it
is sufficient to show that the left-hand side of (\ref{nada}) is
continuous in $t$.

It follows from Theorem 2.18 in~\cite{CFKZ} that for any $z$, $\tilde
{Z}(z,t,r_t(\omega))$ is continuous and has the additivity property
$\Pro_m$-a.e. for on $\LL0,\zeta\LL$. Hence, thanks to (\ref{Zpos})
for $\mathrm{d} z$-a.e. $z$, $\tilde{L}^z_t$ is increasing. One shows then by
monotone convergence that for any positive Borel function~$f$,
$t\rightarrow\int_{\R}f(z)\tilde{L}_t^{z}\di z$ is continuous $\Pro
_m$-a.e. on $\LL0,\zeta\LL$.

For a function $u$ in $\F_{\mathrm{loc}}$, take an nest of closed sets $\{D_k\}
$ and a sequence $(u_k)_{n \in\N}$ of bounded elements of $\F$ such
that $u=u_k$ for q.e. $x\in E$. For any $k\in\N$, let $\tilde
{L}^z_t(u_k)$ be the version $\B(\R)\otimes\B(\R_{+})\otimes\F
^m_{\infty}$-measurable of local time obtained above. Then $\tilde
{L}^z_t:=\tilde{L}^z_t(u_k)$ on $t<\tau_{D_k}$ is a $\B(\R)\otimes
\B(\R_{+})\otimes\F^m_{\infty}$-measurable version of $L^z_t$ and
satisfies the occupation time density formula on $[0,\tau_{D_k}[$, for
any $k\in\N$, so it satisfies it on $[0,\zeta[$.
\end{pf}

\begin{coro}
For any real $a$, $L^a$ is a PCAF and $\Pro_x$-a.e. for q.e. $x\in E$,
the measure in $t$, $\di_t L_t^a$ is carried by the set $\{
s\dvt u(X_{s-})=u(X_s)=a\}$. \label{loctime}
\end{coro}

\begin{pf}
We use $u_k$ and $\{D_k\}$ defined as in the end of the proof of
Proposition~\ref{occuptime}. Since we need to show the assertion of
Cororally~\ref{occuptime} only on $[0,\tau_{D_k}[$, we can assume
that $u$ is a bounded element of $\F$.
It follows from the occupation time density formula and the $\B(\R
)\otimes\B(\R_+)\otimes F^m_{\infty}$-measurability of $\tilde
{L}$, that there exists a subset $R$ of $\R$ of Lebesgue's measure
zero, such that for any $a$ outside of $R$: $\Pro_m$-a.e. $\tilde
{L}_t^a\geq0$ on $\LL0,\zeta\LL$. Consequently, $L^a$ has the same
property. This property holds for any $a\in\R$. Indeed for any real
$a$, take a sequence $(a_n)_{n\in\N}\subset\R\setminus R$ such that
$a_n\downarrow a$. We have: $e(Z^{a_n}-Z^a)=\int1_{\{a<u(x)\leq a_n\}
}\mu_{\langle M^{u,c} \rangle}(\mathrm{d} x)$, which converges to $0$ as $n$
tends to $\infty$ by dominated convergence. Thus, thanks to Theorem
\ref{NakaoFuku} (taking a subsequence if necessary) $\Gamma
(Z^{a_{n}})$ converges to $\Gamma(Z^a)$ uniformly on any finite
interval of $t$, $\Pro_m$-a.e. On the other hand, for $\Pro_m$-a.e.
$w\in\Omega$, $\int_0^t 1_{\{u(X_s)\leq a_{n}\}}\di N^{u,c}_s(\omega
) $ converges to $\int_0^t 1_{\{u(X_s)\leq a\}}\di N^{u,c}_s(\omega)$
for any $t<\zeta(\omega)$. Consequently, we obtain for $\Pro_m$-a.e.
$\omega\in\Omega$, $L_t^a(\omega)\geq0$ for any $t<\zeta(\omega)$.

It follows from Lemma~\ref{lema01280210} that for any real $a$, $L^a$
is a PCAF on $\LL0,\zeta\LL$. By Remark 2.2 in~\cite{CFKZ}, it can
be extended to a PCAF.

Now defining $f(x)=(x-a)^4$ and $h(x)=(x-a)^41_{\{x\leq a\}}$, it
follows from (\ref{5.6.20FOT}) that $fu$ and $hu$ belong to $\F
_{\mathrm{loc}}$. Moreover, we have:
\[
M^{fu,c}_t=4\int_0^t \bigl(u(X_s)-a\bigr)^3\di M_s^{u,c}     \quad
\mbox{and}\quad
M^{hu,c}_t=4\int_0^t \bigl(u(X_s)-a\bigr)^31_{\{u(X_s)\leq a\}}\di M_s^{u,c}
\]
thus, $\langle M^{fu,c},Z^a\rangle=\langle M^{h u,c},M^{u,c}\rangle$,
and from the definition of the stochastic integral (\ref{integral}) we
have that $\Pro_m$-a.e. on $\{t<\zeta\}$
\[
\int_0^t \bigl(u(X_s)-a\bigr)^4 \di\Gamma(Z^a)_s=\int_0^t \bigl(u(X_s)-a\bigr)^41_{\{
u(X_s)\leq a\}} \di\Gamma(M^{u,c})_s.
\]
By Lemmas~\ref{lema01280210} and~\ref{lema01240210}, we finally
obtain: $\int_0^t(u(X_s)-a)^4 \di L_s^a=0$ $\Pro_x$-a.e. for q.e. $x\in
E$.
\end{pf}

\subsection{Integration with respect to local time}\label{secu5.2}

We fix $u$ an element of $ \F$ satisfying (\ref{sup01}) and set:
$l^a_t=\int_0^t1_{\{u(X_{s-})\leq a\}}\di N_s^{u,c}$. Hence, the local
time at $a$ of $u(X)$ satisfies:
\[
L^a=-2\Gamma^a+2l^a.
\]
For any $\omega\in\Omega$ and $t<\zeta(\omega)$, the function
$z\rightarrow l^z_t(\omega)$ is of bounded variation. The application
defined for the elementary functions by
\[
f\rightarrow\sum_{i=0}^{n-1}f_i(l_t^{z_{i+1}}-l_t^{z_i}), \qquad t<\zeta
\]
can hence be extended to the set of locally bounded Borel measurable
functions $f$ from $\R$ into $\R$ as a Lebesgue--Stieljes integral and
we have:
\[
\int_{\R}f(z)\di_z l^z_t=\int_0^t f(u(X_s))\di N_s^{u,c}  ,   \qquad t<\zeta.
\]
Using the stochastic integral with respect to $\Gamma$, the
application defined for the elementary functions by
\[
f\rightarrow\sum_{i=0}^{n-1}f_i(L_t^{z_{i+1}}-L_t^{z_i}), \qquad t<\zeta
\]
can hence be extended to the set of locally bounded Borel measurable
functions $f$ from $\R$ into $\R$ and we have:
\[
-\frac{1}{2}\int_{\R}f(z)\di_z L^z_t=\int_{\R}f(z)\di_z \Gamma
_t^z-\int_0^t f(u(X_s))\di N_s^{u,c},     \qquad t<\zeta.
\]

\section{Multidimensional case}\label{sec6}\label{SMC}

In this section, we need the following notation. For $d\in\N$,
$x=(x^1,\ldots,x^d),y=(y^1,\ldots,y^d)\in\R^d$, we set $x\leq y$ (resp.,
$x<y$) if and only of $ x^i\leq y^i$ (resp., $x^i<y^i$) for each
$i=1,\ldots,d$ and $]x,y]=\{z\in\R^d\dvt x<z\leq y\}$. The vector $\hat
{x}$ is obtained from $x$ by elimination of its coordinate $x^d$, that
is, $\hat{x}=(x^1,\ldots,x^{d-1})$,
$\widehat{]x,y]}=\{z\in\R^{d-1}\dvt
\hat{x}< z \leq\hat{y}\}$.

Let $\varphi$ be a measurable function from $\R^d$ into $\R$. We
define integration of simple functions with respect to $\varphi$ as
follows. For $f$ a simple function, that is, there exists $x,y\in\R
^d$ such that $f(z)=1_{]x,y]}(z)$ for all $z\in\R^d$:
\begin{eqnarray}
\mbox{if }d&=&1\dvt  \qquad \int_{\R}f(z)\di\varphi(z) =\varphi(y)-\varphi
(x),\nonumber\\
\mbox{if }d&>&1\dvt  \qquad \int_{\R^d}f(z)\di\varphi(z) =\int_{\R
^{d-1}}1_{\widehat{]x,y]}}(z)\di\varphi(z,y^d)-\int_{\R
^{d-1}}1_{\widehat{]x,y]}}(z)\di\varphi(z,x^d).\nonumber
\end{eqnarray}
As an example, if there exist functions $h_i$, $1\leq i\leq d$ such
that $\varphi(z)=\prod_{i=1}^dh(z_i)$, then $\int_{\R^d}f(z)\di
\varphi(z)=\prod_{i=1}^d(h_i(y^i)-h_i(x^i))$.

We extend this integration to the elementary functions $f\dvtx \R
^d\rightarrow\R$ (i.e., $f(z)=\sum_{i=1}^na_if_i(z)$ where $f_i$, $1
\leq i\leq n$, are simple functions and $a_i$, $1 \leq i\leq n$, are real
numbers) by setting
\[
\int_{\R^d}f(z)\di\varphi(z) =\sum_{i=1}^da_i\int_{\R
^d}f_i(z)\di\varphi(z).
\]
An elementary function has many representations as linear combination of
simple functions, but as in the Riemann integration theory, the
integral does not depend on the choice of its representation.

Let $u$ be in $\F^d_{\mathrm{loc}}$ where $\F^d_{\mathrm{loc}} = \{ (u^1,u^2,\ldots, u^d)
\dvt u^i \in\F_{\mathrm{loc}} , 1 \leq i \leq d\}$. Let
$\{D_k\}_{k\in\N}$ be a nest of closed set, $\sigma:=\lim
_{k\rightarrow\infty}\sigma_{E\setminus D_k}$ and $(u_k)_{k\in\N}$
a sequence of bounded elements of $\F^d$ such that $u = u_k$ q.e. on $D_k$.

For any $a$ in $\R^d$ and $i$ in $\{1,2,\ldots,d\}$, we define $Z^a(u^i)$
and $\Gamma^a(u^i)$, respectively, in $\M^c_{f\mbox{-}\mathrm{loc}}$ and $ \Nc
_{c,f\mbox{-}\mathrm{loc}}$ by
\begin{eqnarray}
Z_t^a(u^i)
&=&
\everymath{\displaystyle }\left\{
\begin{array}{l@{ \qquad }l}
\int_0^t1_{\{u_k(X_{s-})\leq a\}}\di M_s^{u^i_k,c}
  &\mbox{for } t\leq\sigma_{E\backslash D_k},\\[3pt]
0&\mbox{for } t\geq\sigma,
\end{array}\right.
 \nonumber\\
\Gamma^a(u^i)&=&\Gamma(Z^a(u^i)).\nonumber
\end{eqnarray}
Thanks to the linearity property of $\Gamma$, we have for any
elementary function $f$:
\[
\int_{\R^d}f(z)\di_z\Gamma_t^z(u^i)=\Gamma_t\biggl (\int
_0^tf(u(X_s))\di M^{u^i,c}_s \biggr).
\]
We extend (\ref{norma}) of Section~\ref{SiwrtG} from $d = 1$ to $d
\geq1$, by defining
for $k\in\N$, the norm $\|\cdot\|_{k}$ on the set of measurable functions
$f\dvtx \R^d\rightarrow\R$
\[
\|f\|_{k}:=\sum_{i=1}^d \biggl(\int_E f^2(u_k(x))\mu_{\langle
M^{u^i_k,c} \rangle}(\mathrm{d} x) \biggr)^{1/2}
\]
and we define the set $\Ii$ with the metric $[\cdot,\cdot]$ as in (\ref
{metric}) of Section~\ref{SiwrtG}. The set of elementary functions is
dense in $\Ii$. We have the following version of Lemma~\ref{int}.

\begin{lema}
The applications $f\rightarrow\int_{\R^d}f(z)\di_z\Gamma^z_t(u^i)$
($1 \leq i\leq d$) defined on the set of elementary functions, can be
extended to the set $\Ii$. This extensions, denoted by $\int_{\R
^d}\di_z\Gamma^z(u^i)$, satisfy:
\begin{enumerate}[(ii)]
\item[(i)] $\int_{\R^d}f(z)\di_z\Gamma_t^z(u^i)=\Gamma
(fu*M^{u^i,c})_t\ \forall t\geq0$, $\Pro_x$-a.e. for q.e. $x\in E$.

\item[(ii)] For $(f_n)_{n \in\N}$ sequence of $\Ii$ such that
$[f_n-f]\rightarrow0$, there exists a subsequence $(f_{n_k})_{k \in\N
}$ such that $\int f_{n_k}(z)\di_z\Gamma^z_t(u^i)$ converges
uniformly on any compact of $[0,\zeta)$ ($[0,\infty)$ if $u\in\F
^d$) to $\int f(z)\di_z\Gamma^z_t(u^i)$ for every $1 \leq i\leq d$
$\Pro_x$-a.e. for q.e. $x\in E$.
\end{enumerate}
\label{intpv}
\end{lema}

With can prove a multidimensional version of Lemma~\ref{defZF} with
the same arguments used in its proof. We have the following
multidimensional It\^{o} formula.

\begin{propo}
Let $u$ be an element of $\F^d$ (resp., $\F_{\mathrm{loc}}^d$) and $F\dvtx \R
^d\rightarrow\R$ a continuous function admitting locally bounded
Radon--Nikodym derivatives $f_i =\partial F/\partial x_i$, $1 \leq i
\leq d$, satisfying the following condition for any $1 \leq i\leq d$
and $k\in\N$
\begin{equation}
\lim_{h\rightarrow0}\int_{E}\bigl\{f_i\bigl(u_k(x)+h\bigr)-f_i(u_k(x))\bigr\}^2\mu
_{\langle M^{u^i_k,c} \rangle}(\mathrm{d} x)=0 . \label{condreg}
\end{equation}

Then, $\Pro_x$-a.e. for q.e. $x\in E$, the process $F(u(X_t)),t\in
[0,\infty)$ (resp., $[0,\zeta)$) admits the decomposition
\begin{equation}
F(u(X_t))=F(u(X_0))+M_t(F,u)+Q_t(F,u)+V_t(F,u),
\end{equation}
where $M(F,u)\in\M_{\mathrm{loc}}$, (resp., $\M_{f\mbox{-}\mathrm{loc}}$) $Q(F,u)\in
\Nc_{c,\mathrm{loc}}$ (resp., $\Nc_{c,f\mbox{-}\mathrm{loc}}$) and $V(F,u)$ is a
bounded variation process given by:
\begin{eqnarray}
M_t(F,u)&=&M_t^d(F,u)+\sum_{i=1}^d\int_0^t f_i(u(X_s))\di
M^{u^i,c}_s,\nonumber\\
Q_t(F,u)&=&\sum_{i=1}^d\int_{\R} f_i(z)\di_{z} \Gamma
_t^{z}(u_i)+A_t(F,u),\nonumber\\
V_t(F,u)&=&\sum_{s\leq t}\{F(u(X_s))-F(u(X_{s-}))\}1_{\{
|u(X_s)-u(X_{s-})|\geq1\}}1_{\{s<\xi\}}\nonumber\\
&&{}-F(u(X_{\xi-}))1_{\{t\geq\xi\}}.\nonumber
\end{eqnarray}
\end{propo}

\begin{pf}
As in the proof of Theorem~\ref{itoSMP}, we can assume that $u$ is a
bounded element of $\F$ and each $f_i$ is bounded. For $\phi\dvtx \R
^d\rightarrow\R$ an infinitely differentiable function with compact
support, the function $F_n$ defined by $F_n(z):=\int_{\R
^d}F(z+y/n)\phi(y)\di y$ converges pointwise to $F(z)$. Setting:
$f_{n,i}=\partial F_n/\partial x_i$ we obtain thanks to (\ref{condreg}):
\[
\lim_{n\rightarrow\infty}\int_{E}[f_{n,i}(u(x))-f_i(u(x))]^2 \mu
_{\langle M^{u^i,c} \rangle}(\mathrm{d} x)=0.
\]
The rest of the proof follows step by step the proof of Theorem~\ref{itoSMP}.
\end{pf}

In the case where $E=\R^d$ and $\E^{(c)}$ is given by
\[
\E^{(c)}=\sum_{i,j=1}^d\int_{\R^d}\frac{\partial u}{\partial
x_i}\frac{\partial v}{\partial x_j}a_{ij}(x)\di x,
\]
where for every $(i,j)$, $a_{ij}$ is a bounded measurable function. The
coordinates functions $\pi_i(x)=x_i, 1 \leq i \leq d$, belong to $\F
_{\mathrm{loc}}$ and $M=(M^{\pi_1,c},\ldots,M^{\pi_d,c})$ is a martingale
additive functional with quadratic covariation $\langle M^i,M^j\rangle
_s=\int_0^ta_{ij}(X_s)\di s$, hence,
$\mu_{\langle M^{i,c}\rangle}(\mathrm{d} x)=a_{ii}(x)\di x$, and the
condition (\ref{condreg}) holds for any locally bounded measurable function.

\section*{Acknowledgements}
I would like to thank the referee for his/her suggestions. I would also
like to thank sincerely my thesis advisor Nathalie Eisenbaum for every
helpful discussion that led to improvement of the results in this paper.


\printhistory

\end{document}